\documentclass[12pt]{article}
\usepackage{amsthm}
\usepackage{amsmath}
\usepackage{amssymb}
\usepackage{fullpage}
\usepackage[parfill]{parskip}    
\usepackage{latexsym}
\usepackage[footnotesize,bf]{caption}
\usepackage{url}
\usepackage{graphicx}
\usepackage{gensymb} 

\newcommand\tab[1][1cm]{\hspace*{#1}}  
\newtheorem{thm}{Theorem}[section]
\newtheorem{lem}[thm]{Lemma}
\newtheorem{prop}[thm]{Proposition}
\newtheorem{cor}[thm]{Corollary}
\newtheorem{conj}[thm]{Conjecture}

\begingroup
    \makeatletter
    \@for\theoremstyle:=definition,remark,plain\do{%
        \expandafter\g@addto@macro\csname th@\theoremstyle\endcsname{%
            \addtolength\thm@preskip\parskip
            }%
        }
\endgroup

\title{Conditions for Obtaining Nontrivial Knots from Collections of Vectors}
\author{Joseph Borgatti}
\date{December 13, 2016}   

\begin{document}
\maketitle

\begin{abstract}
We explore under what conditions one can obtain a nontrivial knot, given a collection of $n$ vectors. First, we show how to get a crossing from any 3 vectors equal in magnitude, by arbitrarily picking 2 vectors and identifying the sufficient and necessary criteria for picking a third vector that will guarantee a crossing when the vectors are reordered. We also show that it's always possible for a set of vectors to be reordered to form the unknot, if they sum to $\vec{0}$ when added together.

\tab Our main results are restricted to sets of $n$ vectors that, when reordered appropriately, project to a regular $n$-gon in $\mathbb{R}^2$. We prove that if $n=6$, we cannot form a nontrivial knot with our vectors. The first nontrivial knot possible ($3_1$) is when $n=7$, and the first $4_1$ knot possible is when $n=8$. We prove that if $n\geq7$, we can always reorder the vectors to get a projection of a nontrivial knot, and also provide an algorithm to choose how to reorder the vectors to get such a knot.
\end{abstract}
\section{Introduction}

\begin{figure}[!h]
    \centering
    \includegraphics[scale=0.9]{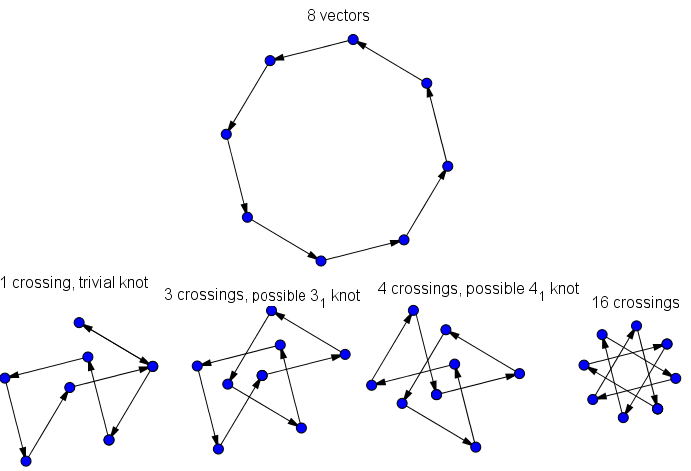}
    \caption{The above collection of 8 vectors can be reordered in many ways, which can result in projections of knots with various levels of complexity, some of which are shown on the bottom. The term ``possible" knot is used, since the $z$-values of the vectors in $\mathbb{R}^3$ must be appropriately chosen.}
    \label{8gons2}
    \end{figure}

	\tab The stick representation of mathematical knots ordinarily are composed of line segments that are as long or short as desired, but cannot bend in any way. In this paper, we instead think of the sticks as collections of vectors, and apply certain limitations on the vectors to see what kind of knots we can get. We  focus on establishing criteria and restrictions for obtaining a projection of a knot with at least 3 crossings from some set of $n$ vectors in $\mathbb{R}^3$, when these vectors are projected onto $\mathbb{R}^2$ and arranged in some reordering. Here, reordering a set of projected vectors means placing them tip-to-tail in a certain order. The reason this topic is interesting is because the same set of vectors can form vastly different knot projections if reordered in different ways, sometimes even by slightly different reorderings. See Figure \ref{8gons2} for an example using a collection of $8$ vectors. One known necessary condition is that any collection of $n$ vectors must sum to $\vec{0}$ when projected into $\mathbb{R}^2$, in order to form a knot. This is because a knot by definition must be a closed loop.
    
	\tab Recall that the three Reidemeister moves, shown in Figure \ref{reide} below, are methods for changing the projection of a knot. These moves, however, do not change the knot itself. Crossings can be called reducible if they can be removed via a series of Reidemeister moves. 
    
    \begin{figure}[!h]
    \centering
    \includegraphics[scale=0.25]{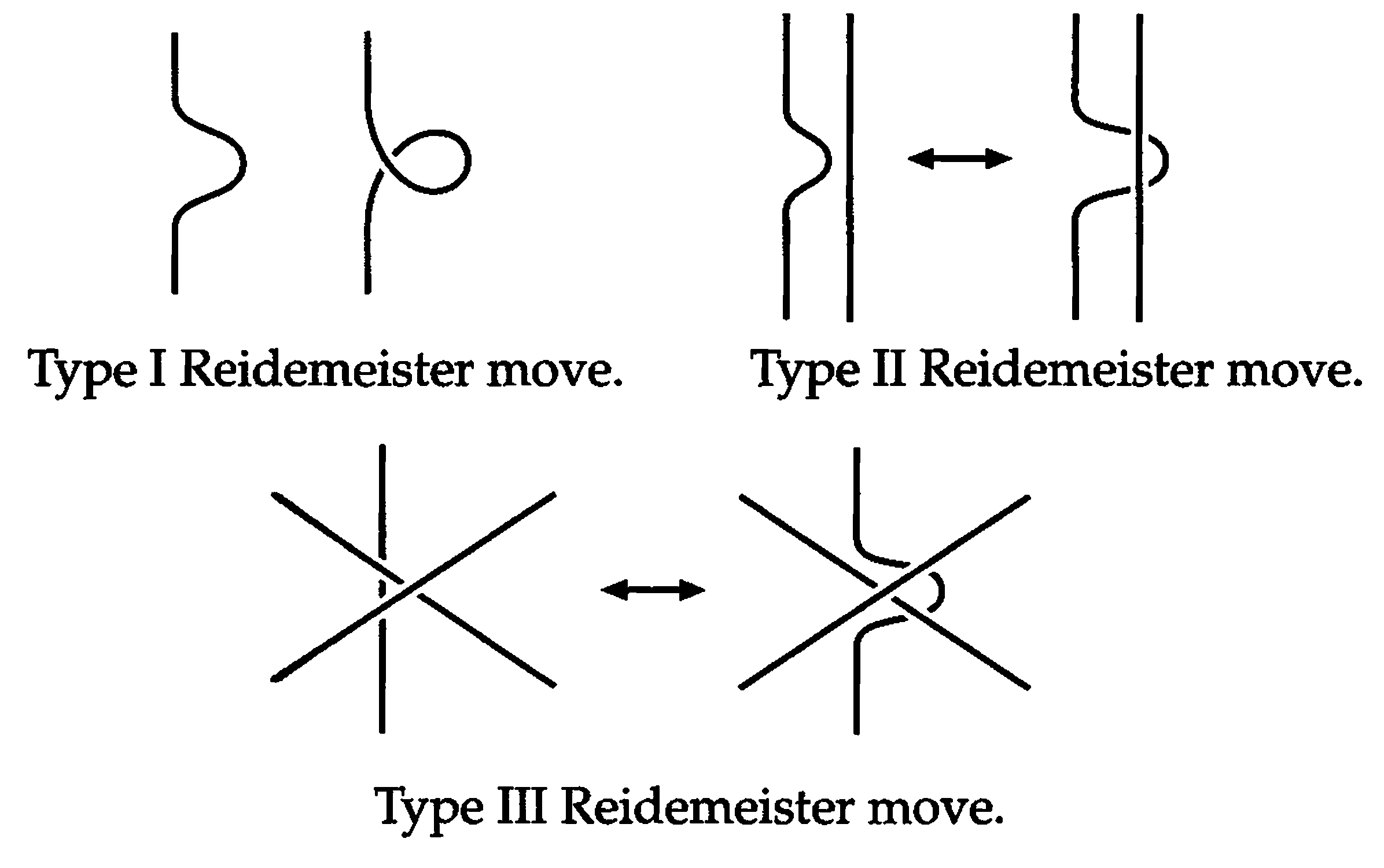}
    \caption{The three Reidemeister moves that can change the projection through ambient isotopy without changing the knot itself. The pictures demonstrate the moves for under-crossings, but they work on over-crossings as well. Figure taken from Adams [2].}
    \label{reide}
\end{figure}
    
    \tab Recall that a knot crossing, sometimes referred to in this paper as a  ``traditional crossing", refers to any point on a knot projection where exactly two strands intersect. The following three knot invariants are used throughout the paper. The crossing number of a knot, $c(K)$, is the minimum number of crossings in any projection of a knot. This means if some knot projection has $n$ crossings, but one or more of the crossings can be removed by a series of Reidemeister moves, then $c(K)<n$ for that knot. The stick number of a knot, $S(K)$, is the minimum number of line segments that can be used to construct a knot. The stick number of the trefoil knot is 6, meaning that if we have less than 6 straight line segments, we cannot make a nontrivial knot. The bridge index of a knot, $br(K)$, is the minimum number of local maximums of a knot from any direction in the xy-plane. The bridge index of a knot is 1 if and only if it is isotopic to the unknot. The values of these invariants for both the unknot, and the prime knots with up to 5 crossings, are shown in \textbf{Table 1} below.

\begin{table}[!h]
\centering
\begin{tabular}{|c|c|c|c|} 
	\hline 
	Knot Type & $c(K)$ & $S(K)$ & $br(K)$\\ 
    \hline
    Unknot & 0 & 3 & 1 \\
    \hline
    $3_1$ & 3 & 6 & 2\\
    \hline
    $4_1$ & 4 & 7 & 2\\
    \hline
    $5_1$ and $5_2$ & 5 & 8 & 2\\
    \hline
\end{tabular}
\caption{A table of knot invariants for the unknot, and for the prime knots with up to 5 crossings. $c(K)$, $S(K)$, and $br(K)$ represent crossing number, stick number, and bridge index (respectively). Values of invariants can be checked at [4].}
\label{invariant table}
\end{table}


	\tab In section 2, we identify some criteria of vectors that will guarantee nontrivial knots when the conditions are met. One result dealt with a regular 6-gon, which led us to investigate what types of knots we could get with higher $n$-sided polygons in section 3. Finally, in section 4 we take a quick look at what types of knots one can get with a triple crossing and a ``traditional crossing" in the knot projection.

\section{Vectors}

In this section, we identify necessary criteria for a projection of a set of $n$ vectors must have in order to make a nontrivial knot in $\mathbb{R}^3$. We start with Lemma \ref{i j}, as in order to make a nontrivial knot we need to make crossings.

\begin{lem} \label{i j} Any collection of vectors must have $i$- and $j$-components of both signs in order to have any crossings.
\end{lem}
This is because by definition, a crossing occurs when we are following the orientation of a knot (or set of vectors) and go back upon a previous strand (or vector).

Corollary \ref{one direction} naturally follows from Lemma \ref{i j}:
\begin{cor} Given a projection of nonzero vectors that sum to $\vec{0}$ in $\mathbb{R}^2$, if only one vector has an $i$-component (or $j$-component) of one sign and the rest have zero or oppositely signed $i$-components (or $j$-component, respectively), then the vectors will form only trivial knots no matter how they are reordered.
\label{one direction}

\begin{proof} Start with any collection of nonzero vectors that, when projected onto $\mathbb{R}^2$, sum to $\vec{0}$. This means that placing them tip-to-tail, a closed loop will start with $\vec{v_1}$ and end at the tail of $\vec{v_1}$. Without loss of generality, let the vector with unique $i$-component (or $j$-component) be $\vec{v_1}$. To get crossings, we would need to go back in a direction towards a previous vector when adding tip-to-tail. However, since we have already used our vector $\vec{v_1}$, any subsequent vectors must go in the opposite $i$ direction of $\vec{v_1}$. This means that any pair of vectors in $\vec{v_2}$...$\vec{v_n}$ will not cross, no matter how they are ordered. 

	\tab It is indeed possible to cross $\vec{v_1}$ by facing vectors towards its $y$-coordinate (or $x$-coordinate, respectively) value. No matter how these crossings look in $\mathbb{R}^3$, either end of the twist could easily be undone with a Reidemeister I move, which could be repeated with the resulting twist, until we are left with the unknot. An example of this is shown in Figure \ref{trivial twist} below. Therefore, no matter how we try to cross a collection of vectors with these constraints, we will always form trivial knots.
\begin{figure}[h!] 
\centering
  \includegraphics[scale=0.6]{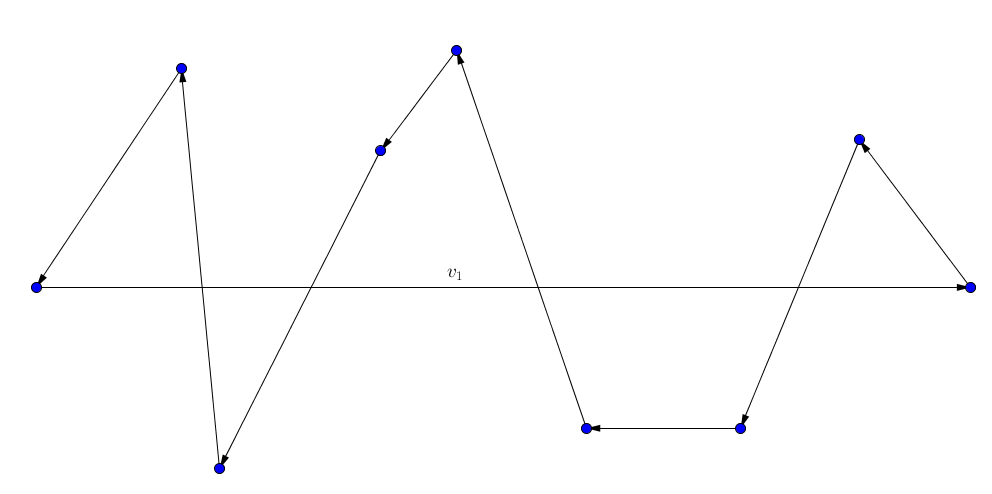}
  \caption{A trivial ``twist" from obtained from a set of vectors, where one vector ($\vec{v_1}$ here) has a uniquely-signed $i$- or $j$-component.}
  \label{trivial twist} 
  \end{figure}
\end{proof}
\end{cor}

	\tab Note that the terms ``$i$- and $j$-components'' referred above do not need to be restricted to a traditional drawing of the xy-plane. Lemma \ref{i j} applies to any orientation or viewpoint for projections of vectors in an xy-plane, as this can be viewed as rotating vectors, which will not affect crossing information. We explore the results of Lemma \ref{one direction} further by investigating what sufficient and necessary conditions three vectors must have in order to make a crossing (even if the crossing is trivial).

\begin{lem} If a collection of three vectors projected onto $\mathbb{R}^2$ have equal length, and are not all in some half-plane when their tails are placed at the origin, then they can be positioned tip-to-tail so that they will either cross or form a closed loop.
\label{vec origin}
\begin{proof}Start with a projection two vectors onto $\mathbb{R}^2$ that have equal length and aren't colinear in $\mathbb{R}^2$. Without loss of generality, call these $\vec{A}$ and $\vec{B}$, and place the tails of these vectors at the point (0,0). Next, draw straight lines extending from the base of each vector. These are depicted by the dotted lines in Figure \ref{3 vec cross 1} below. Each dotted line with its corresponding base vector can be thought of as a boundary for the half of the xy-plane that contains $\vec{A}$ and $\vec{B}$. According to Lemma \ref{i j}, our third vector cannot be in the same half of the xy-plane as $\vec{A}$ and $\vec{B}$, since we would not have $i$- and $j$-components of both signs. This means the third vector, $\vec{C}$, must be exclusively between the two dotted lines. The left image in Figure \ref{3 vec cross 1} shows one option for $\vec{C}$, which results in a crossing when reordered like in the right image.

\begin{figure}[!h]
\centering
\includegraphics[scale=0.9]{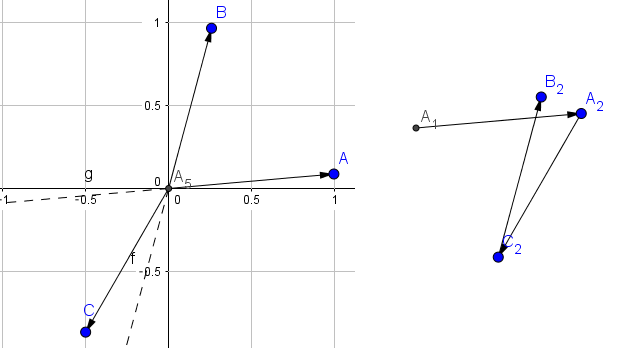}
\caption{On the left, projections of three vectors with tails placed at the origin that can make a crossing if ordered properly. The possible vectors for the third vector $\vec{C}$ are bound between the two dotted lines. On the right, a reordering of these vectors that results in a crossing.}
\label{3 vec cross 1}
\end{figure}

\tab The only instance where the vectors won't cross using this algorithm is when the angle between each pair of vectors is equal. As shown in Figure \ref{3 vec cross 2}, $\vec{C}$ would still touch $\vec{A}$, but would not result in a crossing.
\begin{figure}[!h]
\centering
\includegraphics[scale=0.9]{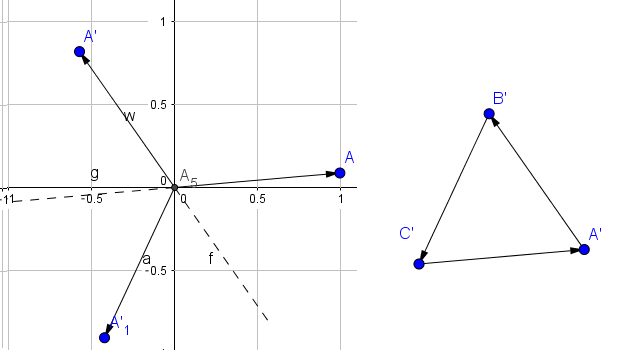}
\caption{This figure is very similar to Figure \ref{3 vec cross 1}, but angle between each pair of the three vectors are equal. This does not give us a crossing, but rather a closed loop that is an equilateral triangle.}
\label{3 vec cross 2}
\end{figure}

\tab Any set of 3 vectors of equal length must cross (or form a closed loop), as without loss of generality there will always be at least one acute angle between either $\vec{B}$ and -$\vec{A}$ (the dotted line opposite to $\vec{A}$) or $\vec{C}$ and -$\vec{A}$. This is true because the angle between $\vec{B}$ and $\vec{C}$ must be less than 180 degrees, since we have limited our option of $\vec{C}$ to be exclusively between the two dotted lines. This means between -$\vec{A}$ and either $\vec{B}$ or $\vec{C}$, there could be either an obtuse angle and an acute angle, or two acute angles. Because when reordering our vectors we can always form an acute interior angle with $\vec{A}$ and our second vector, we can always pick the smallest angle between $\vec{A}$ and either $\vec{B}$ or $\vec{C}$, which will result in the minimal change in y-value. Because our 3 vectors have equal lengths, this means that if $\vec{A}$ is placed on the x-axis, then the third vector will have a greater magnitude in y-value than the preceding vector. Therefore, the third vector must cross (or touch, if equiangular) $\vec{A}$.
\end{proof}
\end{lem}

\tab Next, instead of looking how we can reorder equal-sized vectors to get crossings, in Theorem \ref{unknot guarantee} we look at how one can reorder any set of vectors to get no crossings and guarantee a projection of the unknot.

\begin{thm} Given any set of $n$ vectors in $\mathbb{R}^3$ that sum to $\vec{0}$ when projected onto $\mathbb{R}^2$, it is possible to order them $\vec{v_1}$, $\vec{v_2}$...$\vec{v_n}$ such that the resulting closed curve is the unknot. 
\label{unknot guarantee}


\begin{proof}
	Start with $n$ vectors projected onto $\mathbb{R}^2$ that sum to $\vec{0}$. Imagine the base of all of the vectors being at (0,0) on the xy-plane. Let the polar angle $\theta$ be the counterclockwise angle measured starting from the positive x-axis. Reorder the vectors in terms of ascending value of $\theta$, as shown in Figure \ref{zero sum}. This will guarantee a projection of the unknot in $\mathbb{R}^2$, as the resulting figure will only have one local maximum regardless of orientation. This means the bridge index $br(K)=1$, and since a bridge index of 1 is an invariant only of the unknot, such a reordering of vectors must be the unknot.
    \begin{figure}[h!]
    \centering
    \includegraphics[scale=0.9]{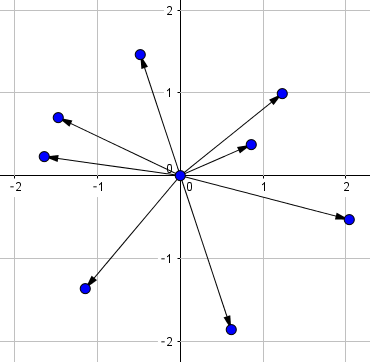}
    \includegraphics[scale=0.9]{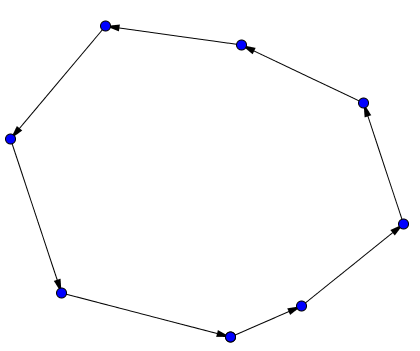}
    \caption{On the left, vectors superimposed on the same base. Reordering these by polar angle $\theta$ guarantees a projection of the unknot (right)}
    \label{zero sum}
    \end{figure}
\end{proof}
\end{thm}

\tab Finally, we go beyond having vectors with equal magnitude (as in Lemma \ref{vec origin}) and add the constraint that the collection of vectors must be able to project to a regular 6-gon in $\mathbb{R}^2$. This means that when the tails of the projected vectors are placed on the origin, the angle between any adjacent pair of vectors will be constant.

\begin{prop} A collection of vectors that form a regular 6-gon in $\mathbb{R}^2$ cannot form a nontrivial knot.
\label{6-gon prop}

\begin{proof}



Start with a closed loop of vectors that form a regular 6-gon in $\mathbb{R}^2$. Like in the algorithm for Lemma \ref{vec origin}, place the base of the 6 vectors at the origin. Because they can form a regular 6-gon, the vectors all have the same length. 

\tab To prove we cannot form a nontrivial knot, we examine the possible reordering of the 6 vectors. Without loss of generality, we can start with vector A. The reordering A(vwxyz), for \{v,w,x,y,z\} $\in$ $\{\vec{B}, \vec{C}, \vec{D}, \vec{E}, \vec{F}\}$, means we start with $\vec{A}$, then add the subsequent vectors tip-to-tail until we return to the tail of $\vec{A}$. The possibilities for the subsequent 5 vectors are encapsulated in the cases below:

\tab \textbf{Case 1: Reordering starts with 3 consecutive vectors.}

\tab If our vector word has 3 consecutive vectors, such as ABC or ABF, then we have used all of the vectors of one $i$ or $j$ direction. By Lemma \ref{one direction}, this means we cannot get a nontrivial knot.

\tab \textbf{Case 2: Reordering starts with the form A\_D.}

\tab Again without loss of generality, start with the reordering ABD. The shape so far resembles a parallelogram without the side parallel to $\vec{B}$. First, consider picking $\vec{C}$. Adding $\vec{E}$ and $\vec{F}$ in either order will not make any crossing. If we instead picked $\vec{E}$ in place of $\vec{C}$, we would form a parallelogram, and the remaining vectors $\vec{C}$ and $\vec{F}$ would be colinear and oppositely directed. If we had picked $\vec{F}$, we would form an equilateral triangle, meaning the remaining vectors would form another equilateral triangle with $\vec{A}$, and thus we would not get any crossings. We do not need to check other options for the second vectors in A\_D, since these would be reflections and rotations of ABD.

\tab \textbf{Case 3: Reordering starts with ACE.}

\tab If the reordering starts with ACE, then this collection of 3 vectors will form an equilateral triangle. This means the remaining vectors, $\vec{B}$, $\vec{D}$, and $\vec{F}$, would form another equilateral triangle no matter how they are ordered, and this would result in 0 crossings.

\tab \textbf{Case 4: Reordering contains a pair of the form AD.}

\tab If the reordering starts with a vector, followed by the vectors that is colinear and oppositely directed to the first one, this does not add any potential to getting crossings. For example, if we start with $\vec{A}$, $\vec{D}$ is the additive inverse of $\vec{A}$, so we will be back at the tail of $\vec{A}$. The remaining 4 vectors must return to $\vec{A}$ in order to form a closed loop, so we now are left with a collection of 4 vectors that start and end at the point of the tail of $\vec{A}$. Because a collection of 4 vectors cannot form a nontrivial knot no matter how they are ordered, this collection of 6 vectors cannot form a nontrivial knot.

\tab We do not need to look at reorderings starting with AE or AF, since these are just isotopic reflections of reorderings starting with AB and AC (respectively). Therefore, no matter how the vectors are reordeded, we cannot form a nontrivial knot from a collection of 6 vectors that project to a regular 6-gon in $\mathbb{R}^2$.
\end{proof}
\end{prop}
\begin{figure} [h!]
\centering
\includegraphics[scale=1.0]{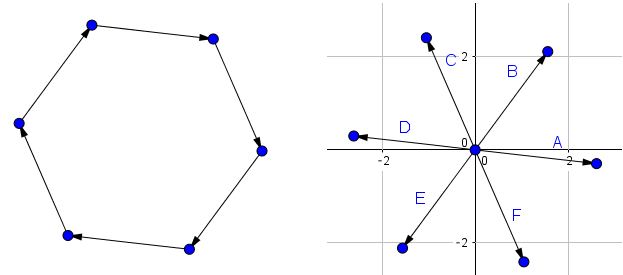}
\caption{Vectors that form a regular 6-gon in $\mathbb{R}^2$, whose bases are all placed on (0,0) on the right.}
\label{6 regular vec}
\end{figure}
Obtaining Proposition \ref{6-gon prop}'s result through using vectors that form a regular 6-gon lead us to investigate what kinds of knots we could get from reorderings of other regular n-gons, for $n\geq7$.

\section{Knots from Regular Polygons}

This section focusing on investigating what types of knots we can get when reordering a collection of vectors that project to a regular n-gon in $\mathbb{R}^2$, for $n\geq7$.

\tab Since in Proposition \ref{6-gon prop} we just proved 6 vectors that form a regular 6-gon in $\mathbb{R}^2$ can only form a trivial knot, the next logical step is to see what knots are possible using a collection of 7 vectors that form a regular 7-gon in $\mathbb{R}^2$. 

\tab Indeed, if we reorder the vectors like in Figure \ref{7 gon trefoil} below, it's possible to get a trefoil knot if we assign crossing information to be alternating (i.e., as we trace around the knot, we never have two consecutive over-crossings or two consecutive under-crossings).

\tab This means that \textbf{7} is the minimum numbers of vectors that form a regular n-gon in $\mathbb{R}^2$ that are needed to get a nontrivial knot, when reordered in some particular embedding, for $n\geq7$.

\begin{figure}[!h]
\centering
\includegraphics[scale=0.75]{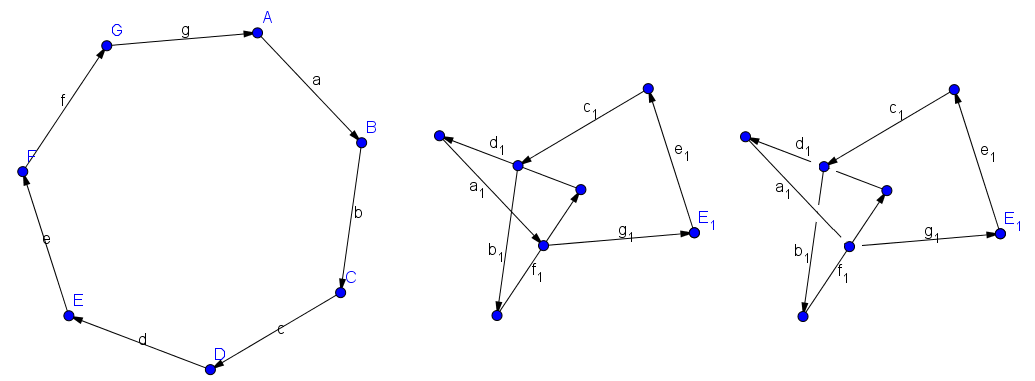}
\caption{On the left, a collection of 7 vectors that project to a regular 7-gon in $\mathbb{R}^2$. In the middle, a reordering of those vectors that gives a graph with 3 crossings. On the right, crossing information is added to make a tentative trefoil knot in $\mathbb{R}^3$.}
\label{7 gon trefoil}
\end{figure}

\tab To prove this the projection on the right of Figure \ref{7 gon trefoil} could be a trefoil in $\mathbb{R}^3$, we imagine point $B$  as being at $(0,0,0)$, and label the $x$ and $y$ values of the remaining points, as shown in Figure \ref{7 vec trefoil heights}. We need to find $z$ values of the intersection points such that the figure in $\mathbb{R}^3$ is not an impossible shape, and thus able to be a trefoil knot. 

\tab One way to do this is to parameterize equations for the lines. First, consider the line segment formed by connecting point $C$ to point $G$, in Figure \ref{7 vec trefoil heights} below. To find how far along the line we must travel to reachthe point (0,0), we write the following equation and solve for ${t_1}$:

\begin{figure} [!h]
\centering
\includegraphics[scale=0.8]{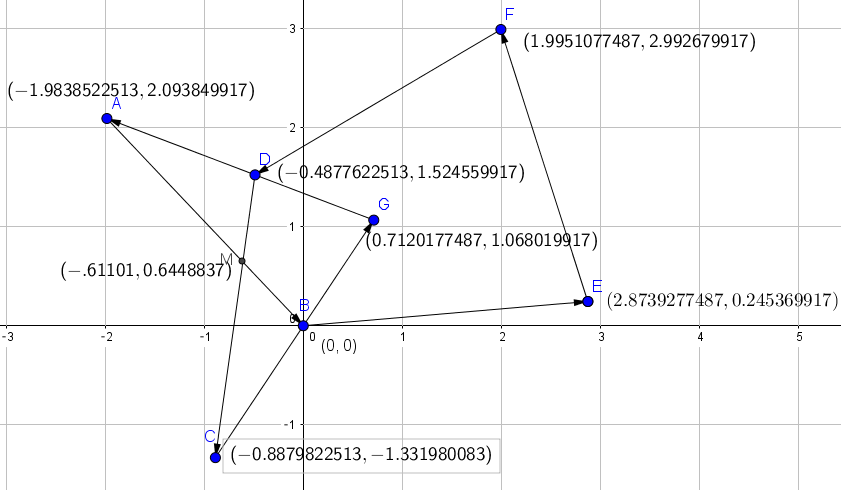}
\caption{Graph from Figure \ref{7 gon trefoil} with added $x$ and $y$ values, if point $B$ was placed at the origin.}
\label{7 vec trefoil heights}
\end{figure}

$(1-t_1)(-.88798,-1.3319) + t_1(0.71201, 1.0680) = (0,0)$

$\implies -.88798 + -.88798t_1 + 0.71201t_1 = 0$ and $-1.3319 + 1.3319t_1 +  1.0680t_1 = 0$

$\implies 1.6t_1 = .88798$

$\implies t_1 = 0.5549889$

This means that at along the line between points $G$ and $C$, at $t_1 = 0.5549889$, 
\[
(1 - t_1)z_C + t_1z_G > 0
\]
where $z_C$ and $z_G$ represent possible $z$-values for points $G$ and $C$, respectively.

We also do this for the segment between $G$ and $A$, and the segment between $A$ and $B$. Using the same process for the other points of intersection, we get
\[
(1-t_2)(z_G) + t_2(z_A) < z_D
\]
\[
(1-t_3)(z_A) > z_M
\]

To find the values of $t_2$ and $t_3$:

$(1-t_2)(0.71201,1.06801) + t_2(-1.98385,2.09384) = (-0.487776, 1.52455)$

$\implies 0.71201 - 0.71201t_2 - 1.98385t_2 = -.487776$

$\implies t_2 = 0.445043715$

and

$(1-t_3)(-1.98385,2.09384) = (-.61103, 0.64488837)$
$\implies t_3 = 0.6919982324$

This gives us a system of three inequalities for possible $z$-values of these points of intersection:

$(0.4450111)z_C + 0.5549889z_G > 0$

$(0.554956285)(z_G) + 0.445043715(z_A) < z_D$

$(0.3080017676)(z_A) > z_M$

One of many solutions to this system of inequalities is $z_C$ = 1, $z_G$ = -0.7. $z_A$ = 7, $z_D$ = 3, and $z_M$ = .1.

This works because $0.05651887 > 0, 2.726836606 < 3$, and $2.156012375 > 0.1$. The $z$-value of the remaining points $F$ and $E$ do not matter.

\tab Because we do, in fact, have a projection of a trefoil knot in $\mathbb{R}^2$, this means that the lowest number of vectors that sum to a regular n-gon in $\mathbb{R}^2$ needed to get a nontrivial knot in $\mathbb{R}^3$ in 7 vectors.

\begin{figure} [h!]
\centering
\includegraphics[scale=0.8]{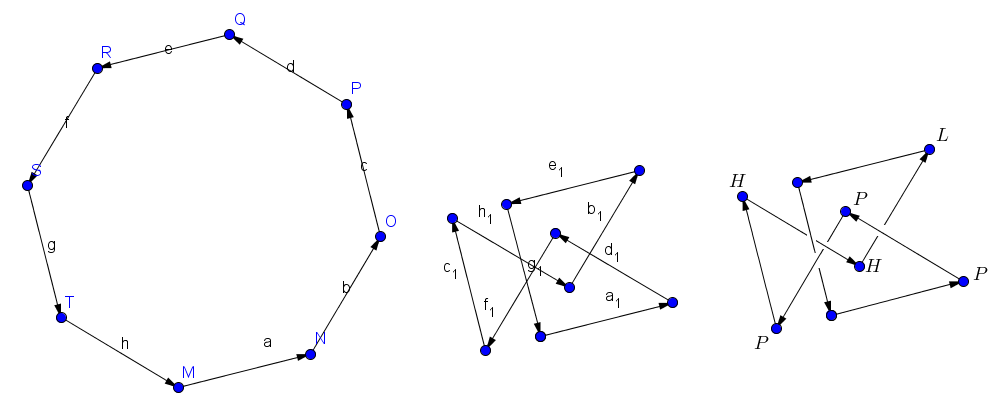}
\caption{On the left, a regular octagon using 8 vectors. In the middle, a reordering of these vectors that form 4 crossings. On the right, crossing information is added that results in a $4_1$ knot in $\mathbb{R}^3$. For the two vertices not labelled, the top can be sufficiently high enough and the bottom can be sufficiently low enough in order for the depicted crossings to be possible in $\mathbb{R}^3$.}
\label{4_1 8 vec}
\end{figure}

\begin{lem} There exists a collection of 8 vectors that form a regular 8gon in $\mathbb{R}^2$, which can be reordered into a projection of the $4_1$ knot.

\begin{proof} Start with a collection of 8 vectors $\vec{v_1}, \vec{v_2} ... \vec{v_8}$ that, when added tip-to-tail, form a closed regular octagon when projected in $\mathbb{R}^2$. These can be reordered in a way that is isotopic to the embedding with 4 crossings, shown in Figure \ref{4_1 8 vec}. Add crossing information so we have an alternating projection with 4 crossings. 

\tab Next, we specify the height of these 8 ``sticks" in $\mathbb{R}^3$. In Figure \ref{4_1 8 vec}, we denote vertices as ``L", ``P", or ``H" to denote that end of the stick being lower than the plane, at the plane, or higher than the plane (respectively) for some flat plane $z$. For example, a vector that starts at ``H" and points towards ``L" has the possibility to undercross a vector at height ``P", if we make ``L" low enough relative to ``P". For the two vertices not labeled with heights, consider a straight upward stick coming out of the page, placed in the triangle formed by the edges $h_1$, $c_1$, and $f_1$. Imagine placing the sufficiently-long stick vertically on the projection, then tilting it so it goes over the edge $h_1$ and under the edge $f_1$. The heights of these two vertices do not matter with respect to the adjacent ``L" and ``P" vertices, since there are no crossings between each respective pair. Therefore, we have formed an alternating knot with 4 crossings in $\mathbb{R}^3$, which is the knot $4_1$.
\end{proof}
\end{lem}

\tab Next, in Conjecture \ref{8 vec counterproof}, we instead see if we can get a projection of the $5_1$ or $5_2$ knot. We obtain a reordering that has 5 crossings in $\mathbb{R}^2$, but trying to assign heights to the vertices of the 8 vectors cannot make a knot more complex than the $4_1$ knot in $\mathbb{R}^3$.

\newpage \begin{conj} A collection of 8 vectors that project to a regular 8-gon in $\mathbb{R}^2$ cannot be reordered to form a $5_1$ knot in $\mathbb{R}^3$, even if a projection with 5 or more crossings is found. 
\label{8 vec counterproof}
\end{conj}
\tab For example, start with a collection of 8 vectors that form a closed regular 8-gon when projected into $\mathbb{R}^2$. These can be reordered to form a closed figure with 5 crossings as shown in Figure \ref{5_1 counterproof} below. 
\begin{figure}[h!]
\centering
\includegraphics[scale=1.0]{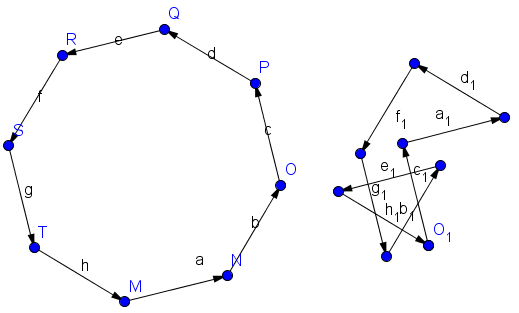}
\caption{A reordering of 8 vectors that gives 5 crossings, but can't be the $5_1$ knot}
\label{5_1 counterproof}
\end{figure}

Notice the three vectors at the top, $a_1$, $f_1$, and $d_1$ are not included in any of the 5 crossings. This means you could replace the path ($a_1$,$d_1$) with just one edge, given by the vector ($a_1 + d_1$) This, however, would result in a figure with 7 sticks. Because the stick number of $5_1$ is 8, we wouldn't have enough sticks. Therefore this reordering of vectors would not have the possibility of being a $5_1$ knot in $\mathbb{R}^3$. 

\tab However, there may be some other reordering that makes 5 crossings and has the possibility of being a $5_1$ knot in $\mathbb{R}^3$, so we cannot say the case for the above example is always true. This lead us to investigate if there is an example where we can get a projection of the $5_1$ with vectors that make a regular n-gon in $\mathbb{R}^2$.

\begin{cor} There exists a set of 8 or more vectors who project to a regular 5-gon in $\mathbb{R}^2$, that can be reordered to create the $5_1$ knot in $\mathbb{R}^3$.
\begin{proof} Start with at least 5 vectors that, when projected in $\mathbb{R}^2$, add to a regular 5-gon. Reorder the vectors tip-to-tail to create the figure below. This is done by starting with some vector, calling it $\vec{v_1}$, and then adding every 3rd vector of the 5gon as the subsequent vector in the reordering. Next, add crossing information to create an alternating projection with 5 crossings. 

\tab However, since $S(5_1) = 8$, we need at least 3 more vectors for our $5_1$ knot to not be an impossible figure in $\mathbb{R}^3$. By adding vectors that are $\vec{0}$ in $\mathbb{R}^2$ yet have verticality in $\mathbb{R}^3$, we will still form a regular 5-gon in $\mathbb{R}^2$. Therefore, at the 3 points that have two labels in Figure \ref{5_1 stick} below, place a straight-up or straight-down vector right on each of the vertices. These are used to necessitate bends for preventing an impossible shape.
When labeling the heights of vertices, 'L' can be as low as needed and 'H' can be as high as needed, since the vectors can have any $k$-component and still project to the same shape in $\mathbb{R}^2$. By using the correct choices of 'L', 'P', and 'H' for height labels, such as in Figure \ref{5_1 stick}, the criteria for a $5_1$ knot in $\mathbb{R}^3$ are satisfied.
\begin{figure}[h!]
\centerline{\includegraphics[scale=0.8]{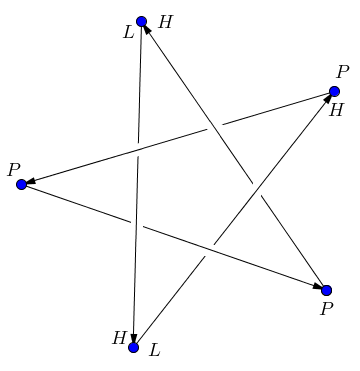}}
\caption{8 vectors that project to a representation of $5_1$ in $\mathbb{R}^2$. 3 of the vectors have the form $< 0, 0, k>$ for $k$ $\in$ $\mathbb{Z}$, and are placed on the vertices that have two height labels. The $k$ values are sufficient enough to facilitate necessary bends to create a $5_1$ knot in $\mathbb{R}^3$.}
\label{5_1 stick}
\end{figure}
\end{proof}
\end{cor}


After looking at what types of knot projections we can obtain $n$=6, 7, or 8 vectors that project to a regular $6$-gon, $7$-gon, or $8$-gon (respectively), we formulate an algorithm that will always guarantee a nontrivial knot projection by reordering $n$ vectors that can project to a regular $n$-gon, for any $n\geq7$.

\begin{figure}
	\centering
    \includegraphics[scale=1.0]{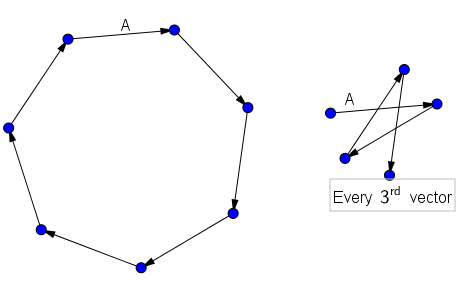}
    \includegraphics[scale=1.0]{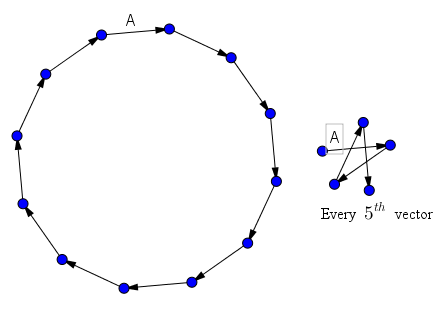}
    \caption{Two examples of choosing which vectors to pick when reordering sets of $n$ vectors that project to a regular $n$-gon in $\mathbb{R}^2$. On the top, we have $n=7$ vectors, so we pick every $\left \lfloor{ \frac{7}{3} }\right \rfloor  + 1 = 3^{rd}$ vector, to get a subset of vectors that have 3 crossings. On the bottom, we have $n=12$ vectors, so we pick every $\left \lfloor{ \frac{12}{3} }\right \rfloor  + 1 = 5^{th}$ vector.}
    \label{n vector generalization}
\end{figure}

\begin{thm} If a set of $n$ vectors project to a regular $n$-gon in $\mathbb{R}^2$, for $n\geq7$, then we can always reorder the vectors to get a projection of a nontrivial knot. 
\begin{proof} To do this, look at the projection of the regular $n$-gon, and start the reordering with any one of the $n$ vectors; call this first vector $\vec{A}$. After we pick $\vec{A}$, we go around the regular $n$-gon and pick every $X^{th}$ vector after the previous vector until we get four vectors, where:

\[
 X = \left \lfloor{ \frac{n}{3} }\right \rfloor  + 1 
 \] 
where $\left \lfloor{ x }\right \rfloor$ represents the floor function evaluated on the real number $x$. 

The motivation for the above relation for $X$ comes from the fact that for very large $n$, the first three $\frac{n}{3}$ vectors would form an equilateral triangle. Adding $1$ to this guarantees that the triangle will be off by even the slightest amount, and would thus would make a crossing instead. The fourth vector that we pick will make two more crossings with two of our previous vectors, resulting in 3 crossings total. Examples of this method are shown in Figure \ref{n vector generalization} later on.

\begin{figure}
	\centering
    \includegraphics[scale=0.8]{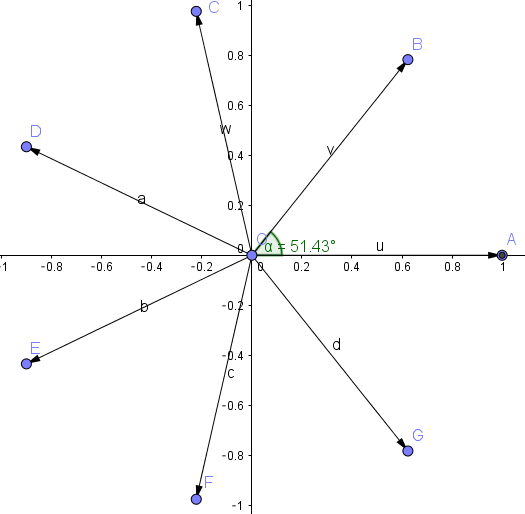}
    \includegraphics[scale=0.8]{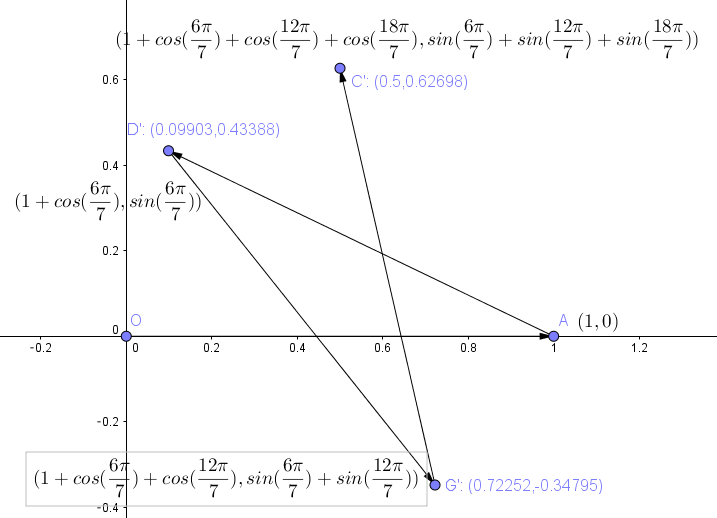}
    \caption{Example of choosing which vectors to pick when reordering a set of $n=7$ vectors of magnitude $1$ that project to a regular $7$-gon in $\mathbb{R}^2$. Since $n=7$, we pick every $\left \lfloor{ \frac{7}{3} }\right \rfloor  + 1 = 3^{rd}$ vector, to get a subset of 4 vectors that have 3 crossings. The order in which we pick the remaining vectors do not matter. On the bottom, coordinates for the tips of each vector are labeled in both decimal form and trigonometric form.}
    \label{n vector generalization 2}
\end{figure}

Without loss of generality, let the length of each vector be $1$, and start the reordering with the vector $<1,0>$. In general, our next vector will be the $X^{th}$ vector after our first one. The tip of the second vector in our reordering of the four vectors will have coordinates of the form:
\[
(1 + \cos(\frac{2\pi}{n}*X), \sin(\frac{2\pi}{n}*X))
\]

In addition, the tip of our third vector will have the following coordinates:

\[
(1 +  \cos(\frac{2\pi}{n}X) +  \cos(2*\frac{2\pi}{n}X) , \sin(\frac{2\pi}{n}X) + \sin(2*\frac{2\pi}{n}X))
\]
And the fourth vector will have the coordinates:
\[
(1 +  \cos(\frac{2\pi}{n}X) +  \cos(2*\frac{2\pi}{n}X) + \cos(3*\frac{2\pi}{n}X) , \sin(\frac{2\pi}{n}X) + \sin(2*\frac{2\pi}{n}X) + \sin(3*\frac{2\pi}{n}X))
\]
To prove that these vectors will always cross, we look at the bounds for $\frac{2\pi}{n}X$. Because $X = \left \lfloor{ \frac{n}{3} }\right \rfloor  + 1$, it follows that:
\[
\frac{n}{3} < X \leq \frac{n}{3} + 1
\]
Multiplying by $\frac{2\pi}{n}$ we get:
\[
\frac{2\pi}{3} < \frac{2\pi}{n}X \leq \frac{2\pi}{3} + \frac{2\pi}{n}
\]
Because we have chosen $n\geq7$, the largest value for $\frac{2\pi}{n}$ is $\frac{2\pi}{7}$. Therefore, our bounds can be sharpened to:
\[
\frac{2\pi}{3} < \frac{2\pi}{n}X \leq \frac{2\pi}{3} + \frac{2\pi}{7} \tab \textrm{or} \tab 120\degree < \frac{360}{n}X\degree \leq 171.4285714 \degree
\]
Let $\phi =  \frac{2\pi}{n}X$. For all $n$, the $y$-coordinate for the tip of our second vector are within the bounds $\sin(\frac{2\pi}{3})$ and $\sin(\frac{2\pi}{3} + \frac{2\pi}{7})$. This is approximately:
\[
0.149042 < \sin(\phi) < \frac{\sqrt[]{3}}{2}
\]
Next, we look at the $y$-coordinate of the tip of the third vector to see if it is always less than $0$, since the $y$-coordinate for any point on the first vector is $0$. We also check the $x$-coordinate, as even if the $y$-coordinate was negative, a negative $x$-coordinate could mean the third vector never actually touches the first vector.
\[
-.369009 < \sin(\phi) + \sin(2\phi) < 0
\]
\[
0 < 1 + \cos(\phi) + \cos(2\phi) < 0.966742
\]
Because for $\frac{2\pi}{3} < \phi \leq \frac{2\pi}{3} + \frac{2\pi}{7}$ the $y$-coordinate of the tip of the third vector is always negative and the $x$-coordinate is always positive, the third vector must cross the first. 

Finally, we look at the tip of the fourth vector and see if its $y$-coordinate it always above the point on the second vector at the same $x$-value. To do this, we first write an equation of a line for the second vector, where $(x_1,y_1)$ is (1,0) and $(x_2,y_2)$ is the tip of the second vector.
\[
m = \frac{y_2 - y_1}{x_2-x_1} = \frac{\sin(\phi)}{(1+\cos(\phi))-1} = \tan(\phi)
\]
Using point-slope form:
\[
y - \sin(\phi) = \tan(\phi) (x - ((1+\cos(\phi))
\]
\[
 => y = x*\tan(\phi) - \tan(\phi) - \tan(\phi)\cos(\phi) + \sin(\phi)
\]
\[
 => y = x\tan(\phi) - \tan(\phi)
\]
Comparing the $y$-value of the tip of the fourth vector and the $y$-value of the point on the second vector at $x = 1+\cos(\phi)+\cos(2\phi)+\cos(3\phi)$ gives us
\[
(\sin(\phi) + \sin(2\phi) + \sin(3\phi)) - (1+\cos(\phi)+\cos(2\phi)+\cos(3\phi)\tan(\phi) - \tan(\phi))
\]
which, after graphing, is always greater than zero. Therefore, the tip of the fourth vector will always be above the corresponding point on the second vector. Because of this, and the fact that the tip of the third vector is always negative. By the Intermediate Value Theorem, this means there exists some point $c$ on the fourth vector, such that the $y$-value of $c$ is between the tip of the third vector and the tip of the fourth, which implies the fourth vector and second vector must intersect. This is also true for the first and fourth vectors. 

Therefore as we have proven there will always be 3 points of intersection among our 4 vectors, this means that we can get a projection of a nontrivial knot by choosing the reordering described above.

The remaining $n - 4$ vectors can be added to our subset of 4 vectors in any order as we do not care about getting a knot more complex than the $3_1$ knot. The simplest way to reorder the remaining vectors would be in terms of ascending $\theta$, where $\theta$ is the polar angle measured from the positive $x$-axis when the tails of the vectors are placed at the origin of the $xy$-plane. The resulting projection will always be isotopic to a projection of the $3_1$ knot. Additionally, if the heights of the vectors in $\mathbb{R}^3$ have sufficiently high and/or low $z$ values, then we can even get a knot isotopic to the $3_1$ knot in $\mathbb{R}^3$.
\end{proof}
\end{thm}

\newpage \tab When looking at what kind of knots we could get with vectors that form a regular 7-gon in $\mathbb{R}^2$, many reorderings formed strange crossings like in Figure \ref{7 vec} below:
\begin{figure} [h!]
\centering
\includegraphics[scale=0.8]{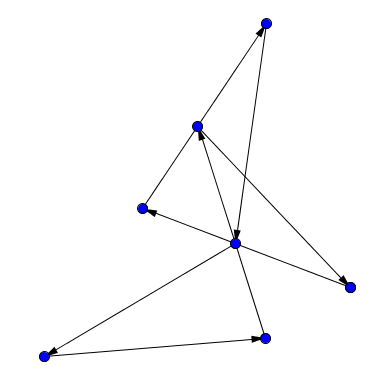}
\caption{Knot projection with a 'triple crossing' obtained by reordering a projection 7 vectors of a regular 7-gon in $\mathbb{R}^2$}
\label{7 vec}
\end{figure}

\tab This leads to Section 4, where we briefly investigate how to interpret crossings of these types.
\section{Triple Crossings}

A \textbf{Triple Crossing} of a knot projection has three strands intersecting at a point, as opposed to traditional crossings which are only composed of two strands. To define what this means for the knot in $\mathbb{R}^3$, one must assign heights to the three strands. Once the heights of the strands are labeled, one can view the triple crossing from a different angle, and interpret it as 3 separate traditional crossings. There are 6 ways we can label the set of strands; one of which is shown on the right of Figure \ref{Triple} below.
\begin{figure}[!h]
\centering
\includegraphics[scale=0.9]{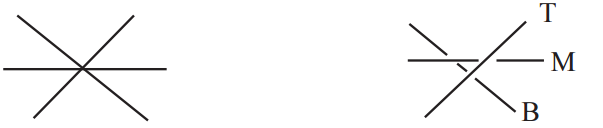}
\caption{On the left, a top-down view of a triple crossing before adding crossing information. On the right, the strands are labeled 'T', 'M', and 'B' for top, middle, and bottom strands (respectively). Notice that the same triple crossing can be viewed as three regular crossings when viewed from a different angle. Figure taken from Adams [1].}
\label{Triple} 
\end{figure}

\newpage \begin{thm} A knot projection with a triple crossing and one other traditional crossing can only be the unknot or trefoil knot.
\begin{proof} Start with a triple crossing like on the left side of Figure \ref{Triple} above, and some other traditional crossing not yet connected to the triple crossing. 
Label the 6 ends of the strands in the triple crossing as $(1,2,3,4,5,6)$ and the 4 ends of the strands in the traditional crossing as $(a,b,c,d)$. Note that any element of $\{a,b,c,d\}$ cannot connect to another element from the same set, as this would either make a link or a Reidemeister I move. This means that each of $\{a,b,c,d\}$ must connect to a unique element of $\{1,2,3,4,5,6\}$. Therefore, one end of a strand in the triple crossing must connect with an end of another, distinct strand of the triple crossing. This greatly limits the possible ways we could connect the strands of the triple crossing and traditional crossing.

\tab Start by connecting the ends of the strands in the triple crossing labeled ``5" and ``6" in Figure \ref{triple cross} below. We can pick these without loss of generality because of rotational symmetry. Now, we must connect each of $\{1,2,3,4\}$ to any unconnected element of $\{a,b,c,d\}$. Again without loss of generality, start with ``1" and connect it to ``a". For each next pair to connect, pick the next clockwise strand of the triple crossing and the next counterclockwise strand of the traditional crossing. This will always ensure that we do not ``block" an end of some strand, which would result in requiring an additional crossing in order to close the loop later on. In our example, we would connect ``2" and ``c", ``3" and ``d", and ``4" and ``b". 


\tab Now that we have a knot projection with a triple crossing and traditional crossing, we must add crossing information. However, we don't have to check all of these cases. No matter how we label the strands of the triple crossing, the top strand will always be part of 2 over-crossings and the bottom strand will always be part of 2 under-crossings. Without creating more crossing, performing a series of Reidemeister moves on the projection will result in a projection isotopic to the unknot or trefoil knot. Examples of each are shown in Figure \ref{triple cross} below. Therefore, the only knots we can get from a projection with one triple crossing and one traditional crossing are the unknot and trefoil knot.
\end{proof}

\begin{figure}[!h]
	\centering
    \includegraphics[scale=1.0]{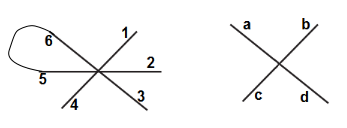}
    
    \includegraphics[scale=0.8]{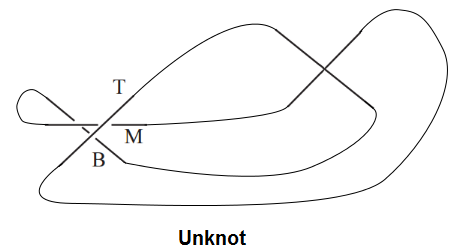}
    \includegraphics[scale=0.8]{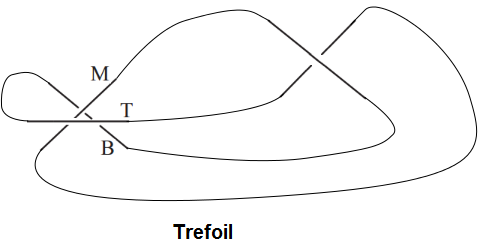}
    \caption{On the top, a triple crossing and a traditional crossing with the ends of each strand labeled. To make a knot projection with 4 crossings, we must connect each of $\{1,2,3,4\}$ to any unconnected element of $\{a,b,c,d\}$. On the bottom left, crossing information is added to form an uknot (crossing information for the traditional crossing does not matter). On the bottom right, crossing information is added but this time we get a trefoil. The traditional crossing for the trefoil must have the depicted crossing information, since the other option would make three consecutive over-crossings which would allow for a Type III Reidemeister move that will untie the projection.}
    \label{triple cross}
\end{figure}
\end{thm}



\newpage
\section{Acknowledgments \& Future Work  }
A main result of the paper involved the creation of an algorithm that will always get a projection of a nontrivial knot, given a collection of $n\geq7$ vectors that project to a regular $n$-gon in $\mathbb{R}^2$. An excellent idea for future work would be to design algorithms for getting even more complex knots, and/or algorithms that do not require the collection of vectors to project to a regular $n$-gon in $\mathbb{R}^2$.

When exploring and reordering vectors that project to a regular $n$-gon in $\mathbb{R}^2$, this paper only dealt with collections of up to 8 vectors, as the term of research was only a semester long. It would be interesting to investigate what types of knots one can obtain from reordering collections of 9 or more of these kinds of vectors. Perhaps with the use of computer processing, one would be able to find a relation between the number of sides of the regular $n$-gon and the maximum crossing number possible for a knot obtained from reordering the vectors. For example, in our paper we proved that the first knot with $c(K)=3$ can be obtained with $n=7$ vectors, and the first knot with $c(K)=4$ can be obtained with $n=8$ vectors. Can we get a $5_1$ or $5_2$ knot with $n=9$ vectors? $n=10$ vectors?

This paper was completed as part of an undergraduate Honors Capstone project at Merrimack College, with supervision and guidance by Dr. Dana Rowland of the Mathematics Department. Figures not sourced were made using the programs GeoGebra and Inkscape.
 
\section{References}
[1] C. Adams, \textit{Triple Crossing Number of Knots and Links},
J. Knot Theory Ramifications, 22 (2013), p. 1350006.

[2] C. Adams, \textit{The Knot Book: An Elementary Introduction to the Mathematical Theory of Knots}, American Mathematical Society, Providence, RI, 2004.

[3] J. H. Conway and C. McA. Gordon,
\textit{Knots and Links in Spatial Graphs},
\newblock J. Graph Theory, {7(4)} (1983), 445-453.

[4] E. W. Weisstein, \textit{Knot Invariants}, Wolfram MathWorld -- A Wolfram Web Resource, 
\newblock mathworld.wolfram.com/topics/KnotInvariants.html.



\end{document}